\documentclass[12pt]{amsart}
 \usepackage{mathtools}
\mathtoolsset{showonlyrefs,showmanualtags}

\usepackage[backrefs]{amsrefs}

 \usepackage[top=1.5in, bottom=1.5in, left=1.25in, right=1.25in]	{geometry}
\usepackage[all]{xy}
\usepackage{amsmath}
\usepackage{amssymb}
\usepackage{amsthm}
\usepackage{amscd}

\newcommand{\RR}{\mathbb  R}

\numberwithin{equation}{section}

\newtheorem{theorem}[equation]{Theorem}

\newtheorem{proposition}[equation]{Proposition}

\newtheorem{lemma}[equation]{Lemma}

\usepackage{hyperref} %,hypertexnames=false,colorlinks,[pagebackref]
\hypersetup{
    colorlinks=true,       % false: boxed links; true: colored links
    linkcolor=blue,          % color of internal links
    citecolor=magenta,        % color of links to bibliography
    filecolor=magenta,      % color of file links
    urlcolor=cyan           % color of external links
}

\begin{document}
\title{An Approach To Endpoint Problems in Oscillatory Singular Integrals}
%Enter your title between curly braces
\author{Alex Iosevich}
\address{Department of Mathematics, University of Rochester, Rochester, NY}
\email{iosevich@gmail.com}
\author{Ben Krause}
\address{
Department of Mathematics,
University of Bristol, Bristol BS8 1UG}
\email{ben.krause@bristol.ac.uk}
\author{Hamed Mousavi}
\address{Department of Mathematics, University of Bristol, BS8 1QU, UK}
\email{gj23799@bristol.ac.uk}

\date{\today}

\maketitle

\begin{abstract}
In this note we provide a quick proof that maximal truncations of oscillatory singular integrals are bounded from $L^1(\mathbb{R})$ to $L^{1,\infty}(\mathbb{R})$. The methods we use are entirely elementary, and rely only on pigeonholing and stationary phase considerations.
    \end{abstract}

\section{Introduction}
Our topic is maximal truncations of oscillatory singular integrals, namely understanding operators of the form
\begin{align}\label{e:trunc} \sup_r \big|\int_{|y| > r} e(P(x,y)) f(x-y) K(y) \ dy \big|, \; \; \; e(t) := e^{2\pi i t},\end{align}
where $P(x,y) \in \mathbb{R}[x,y]$, and $K : \mathbb{R}^d \smallsetminus \{0\} \to \mathbb{C}$ is a Calder\'{o}n-Zygmund kernel. This field was initiated by Ricci-Stein \cite{RS1, RS2}, with endpoint issues addressed by one of us, in collaboration with Lacey, \cite{KL2}, see also \cite{KL1,KLW}, building off deep work of Christ and Chanillo \cite{CC}. $L^p(\mathbb{R}^d)$ bounds for \eqref{e:trunc} for $1 < p  <\infty$ are relatively straightforward to establish, typically proven by interpolating $L^2(\mathbb{R}^d)$-estimates, which in turn arise from stationary phase considerations; we focus on endpoint issue in the one-dimensional setting, namely developing the $L^1(\mathbb{R})$ theory. In the interest of keeping this note as self-contained as possible, we do not attempt to maximize generality, but rather focus on the case where $P \in \mathbb{R}[\cdot]$ is a polynomial of $1$ variable, and $K(t) = \text{p.v. } \frac{1}{t}$ is the Hilbert transform. The full generality, in the one-dimensional setting, can be recovered by arguing as in \cite[Lemma 3.1]{KL2}. The major feature of this note is that our methods are entirely elementary, relying only on the pigeon-hole principle and standard stationary phase estimates, see e.g.\ \cite[Appendix B]{BOOK}. We prove the following result.

\begin{theorem}\label{t:main0}
There exists an absolute constant, $0 < C_d < \infty$, so that for any polynomial $P \in \mathbb{R}[\cdot]$ of degree $\leq d$,
\[ \| \sup_{r} \big| \int_{|t| > r} e(P(t)) f(x-t) \ \frac{dt}{t} \big| \|_{L^{1,\infty}(\mathbb{R})} \leq C_d \| f \|_{L^1(\mathbb{R})}. \] 
\end{theorem}

Some preliminary remarks before beginning the proof proper: first, by modulating appropriately, we may assume that all polynomials in question vanish to at least degree $2$ at the origin; second, by dilation invariance, we may assume that
\[ \| P \| := \sum_{j} |a_j| = 1, \; \; \; P(t) = \sum_j a_j t^j;\]
third, if we let $|\psi| \lesssim 1_{|t| \approx 1}$ be an odd compactly supported Schwartz function that resolves the singularity $\frac{1}{t}$ in that
\[ \sum_i \psi_i(t) := \sum_i 2^{-i} \psi(2^{-i}t) = \frac{1}{t}, \ t\neq 0,\]
it suffices to consider 
\[ \sup_{k \geq 1} \big|\sum_{i = 1}^k \int e(P(t)) f(x-t) \psi_i(t) \ dt \big|,\]
by an argument with the Hardy-Littlewood Maximal function and the weak-type $(1,1)$ boundedness of the maximally truncated Hilbert transform; fourth, if we let
\[ E := \{ k \geq 1 : 2^{-100 d} \leq \frac{|a_j| 2^{kj}}{|a_l| 2^{kl}} \leq 2^{100d} \text{ for some } j \neq k \leq d \} \]
then by sparsification and the estimate $|E| \lesssim d^3 = O_d(1)$, we need only consider
\[ \sup_{k \geq 1} \big|\sum_{1 \leq i \leq k, \ i \notin E, \ i = m \mod 10d} \int e(P(t)) f(x-t) \psi_i(t) \ dt \big| \]
 as we can simply address the scales in $E$ using convexity; in particular, for all $k \notin E$, we may bound
\[ |\int e(P(t)) \psi_k(t) \ dt| \lesssim 2^{-10k}, \]
as for these scales $P(t)$ satisfies the same differential inequalities (up to degree $10$, say) as do pure monomials.

Thus, our task reduces to proving the following proposition.
\begin{proposition}\label{p:main0}
There exists an absolute constant, $0 < C_d < \infty$, so that for any polynomial $P \in \mathbb{R}[\cdot]$ of degree $\leq d$ which vanishes to degree $\geq 2$ at $0$ with $\| P \| = 1$, for each $m \leq 10 d$
\begin{align}\label{e:goal} &\| \sup_{k \geq 1} \big| \sum_{1 \leq i \leq k, \ i \notin E, \ i \equiv m \mod 10d} \int e(P(t)) f(x-t) \psi_i(t) \ dt \big| \|_{L^{1,\infty}(\mathbb{R})} \\
& \qquad \qquad \qquad \qquad \qquad \leq C_d \| f \|_{L^1(\mathbb{R})}. \end{align}
\end{proposition}
In what follows, in the interest of notational ease, we will implicitly restrict all scales considered to be those considered in \eqref{e:goal}, for an arbitrary but fixed $m \leq 10d$.

With these preliminary reductions in mind, we turn to the proof.

\subsection{Notation}
Here and throughout we use $e(t) := e^{2\pi i t}$; $M_{\text{HL}}$ denotes the Hardy-Littlewood maximal function, namely
\[ M_{\text{HL}}f := \sup_{r > 0} \frac{1}{2r} \int_{|t| \leq r} |f(x-t)| \ dt.\]

%This can be accomplished by setting
%\[ \psi(t) = \frac{ \eta(t) - \eta(2t) }{t},\]
%where $1_{[-1/2,1/2]} \leq \eta \leq 1_{[-1,1]}$ is an even Schwartz function. 

For intervals $I$, we will let $f_I$ denote $f \cdot \mathbf{1}_I$, and let $CI$ denote interval concentric with $I$ with $C$ times the length.

We will make use of the modified Vinogradov notation. We use $X \lesssim Y$, or $Y \gtrsim X$ to denote the estimate $X \leq CY$ for an absolute constant $C$, which may vary from line to line. If we need $C$ to depend on a parameter,
we shall indicate this by subscripts, thus for instance $X \lesssim_p Y$ denotes
the estimate $X \leq C_p Y$ for some $C_p$ depending on $p$. We use $X \approx Y$ as
shorthand for $Y \lesssim X \lesssim Y $.

%\begin{lemma}
%    Suppose $P(t) = \sum_{1 \leq \alpha \leq d} \lambda_\alpha t^\alpha \in \mathbb{R}[t]$. Then there is a set $K \subset \mathbb{Z}$ of size $O(d^2 \log A)$ so that for all $k \notin K$
%\[ \max_{\pm} \big( \frac{|\lambda_\alpha| 2^{k\alpha}}{|\lambda_\beta| 2^{k \beta}} \big)^{\pm} %\gg A \]
%for all $\alpha \neq \beta$.
%    \end{lemma}

%Henceforth, every interval we consider will have length $2^k$ for $k \notin K$.

%Suppose that $k \notin K$. Then
%\begin{align}
%    |\int e(P(v+t) - P(t)) \psi_k(t) \ dt| \lesssim \mathbf{1}_{|t| \lesssim 1}(v) + 2^{-k} \mathbf{1}_{|t| \lesssim 2^k}(v)
%\end{align}

%To see this, we may assume that $|v| \gg 1$.

%Since $k \notin K$, there exists a monomial $\lambda_\alpha$ so that
%\[ |\lambda_\alpha| 2^{k \alpha} \gg A \cdot \max_{\beta \neq \alpha} |\lambda_\beta| 2^{k \beta} \]
%so that on the support of $\psi_k$,
%\[ P(v+t) - P(t) = \lambda_\alpha \big( (v+t)^\alpha - t^\alpha \big) + \Omega(v,t) \]
%where
%\[ \Omega(v,t) \in \mathbb{R}[v,t]\]
%has a very small coefficient norm relative to scale
%\[ \sup_{|v| \lesssim 2^k} \| \Omega(v,2^k \cdot ) \| \ll A^{-1} |\lambda_\alpha| 2^{k\alpha}. \]
%Since $P$ is a polynomial, we may assume that its derivative is monotonic, possibly after smoothly partitioning the support of $\psi_k$ into $O(d)$ many sub-intervals.

%\[ |( P(v+t) - P(t) )'| \approx  |\lambda_\alpha| |\alpha| |v| 2^{k(\alpha - 2)} \Omega(v,t) \]

\section{The Proof of Proposition \ref{p:main0}}\label{s:app}
We begin by collecting intervals according to their local density
\[ \mathcal{Q}_t(k):= \{ |I| = 2^k, I \not \subset X : \int_I f \approx 2^{-t}|I| \}, \ X := \{ M_{\text{HL}} f \gtrsim 1\}\]
and set 
\[ \mathcal{Q}_t := \bigcup_{k \geq 1} \mathcal{Q}_t(k) \subset \{ M_{\text{HL}} f \gtrsim 2^{-t}\}.\]

It suffices to show that
\[ |\{ \sup_{k \geq 1} \big| \sum_{1 \leq i \leq k} \sum_{I \in \mathcal{Q}_t(i)} T_I f \big| \gtrsim t^{-C} \}| \lesssim t^{-C} \|f\|_{L^1(\mathbb{R})},\]
where here
\[ T_I f(x) := \int f_I(x-t) e(P(t)) \psi_k(t) \ dt, \ |I| = 2^k.\]

Using the pointwise inequality, $|T_I f| \lesssim 2^{-t} \mathbf{1}_{CI}$, we may discard the first $2^{t/2}$ scales (say), and aim to show that

\begin{equation}\label{goal1}
 |\{ \sup_{k \geq 2^{t/2} } \big| \sum_{i=2^{t/2} }^k \sum_{I \in \mathcal{Q}_t(k)} T_I f \big| \gtrsim t^{-C} \}| \lesssim t^{-C} \|f\|_{L^1(\mathbb{R})}.
\end{equation}

By convexity,
\[ \|T_I f\|_{L^1(\mathbb{R})} \lesssim \|f_I\|_{L^1(\mathbb{R})},\]
but unless the $L^1(\mathbb{R})$-mass of $f_I$ concentrates on small sets, we can improve this estimate by the principle of stationary phase,
\begin{equation}\label{est}
\|T_I f\|_{L^1(\mathbb{R})}^2 \lesssim |I| \| T_I f\|_{L^2(\RR)}^2 = \int \int f_I(x) \cdot K_I(x,y) \cdot f_I(y)  \ dx dy,
\end{equation}
where
\[ K_I(x,y) \lesssim \mathbf{1}_{|x-y| \lesssim 1} + |I|^{-1}.\]
In particular, decomposing $I$ into a disjoint union of dyadic intervals
\[ A_I := \{ J \subset I \text{ dyadic}, \  |J| = 1 \};\]
the point of \eqref{est} is that we may dominate
\[ \| T_I f\|_{L^1(\mathbb{R})} \lesssim \sup_{J \in A_I} \| f_J \|^{1/2}_{L^1(\mathbb{R})} \cdot \| f_I \|^{1/2}_{L^1(\mathbb{R})}.\]
Indeed, if we let
\[ A := \sup_{J \in A_I} \| f_J \|_{L^1(\mathbb{R})},\]
then we may bound the left side of \eqref{est} by
\begin{align}
    A \| f \|_{L^1(\mathbb{R})} + \int \frac{1}{|I|} \sum_{|J| = 1, \ J \subset I} f_J(x) \cdot f_I(y) \ dxdy  \lesssim A \|f_I \|_{L^1(\mathbb{R})}. 
\end{align} 

With this observation in mind, we may take an $\ell^1$ sum over all ``light" intervals $\{ J \in A_I : I \}$ so that 
\[  \| f_J \|_{L^1(\mathbb{R})} \lesssim \left( \log_2 |I| \right)^{-C} 2^{-t} |I|,\]
and seek to prove that 
\begin{equation}\label{goal2}
|\{ \sup_{k \geq 2^{t/2} } \big| \sum_{i=2^{t/2} }^k \sum_{I \in \mathcal{Q}_t(k)} T_I b_I \big| \gtrsim t^{-C} \}| \lesssim t^{-C} \| f\|_{L^1(\mathbb{R})},
\end{equation}
where $b_I := \sum_{J \in A_I \text{ heavy}} f_J$, and heavy intervals are simply intervals that are not light.

\subsection{Pigeonholing and A Simple Case}
To establish \eqref{goal2}, we will need to make pigeonholing arguments, and to first dispose of certain ``harmless" terms.

For each $I \in \mathcal{Q}_t$, collect
\[ A_I^l := \{ J \in A_I \text{ heavy} : \int_J f \approx 2^{-l-t} |I| \}.\]
Note that $A_I^l = \emptyset$ unless $2^{cl} \leq \log_2 |I|$. Define
\[ b_I^l := \sum_{J \in A_I^l} b_J, \; \; \;  b_I^{l,r} := b_I^l \cdot \mathbf{1}_{ \| b_I^l\|_{L^1(\mathbb{R})} \approx 2^{-r-t}|I| },\]
so that we necessarily have $r \leq l$.

It is easy to see that the supports of $\{ b_I^l : I \in \mathcal{Q}_t \}$ are disjoint: if $J \in A_I^l \cap A_{I'}^l$, with $I \subsetneq I'$, this would force
\[ \int_J f \approx 2^{-l-t} |I| \text{ and } \int_J f \approx 2^{-l-t} |I'|,\]
which cannot be, since $|I| < |I'| \Rightarrow |I| \leq 2^{-10 d} |I'|$ by our initial sparsification.

Next, set
\[  k(t,l) := \lfloor \max\{ 2^{t/2}, 2^{cl} \} \rfloor.\]
Moving forward, this will represent our minimal scale; all intervals considered will be of length at least $2^{k(t,l)}$.

Our goal will be to establish the following estimate
\begin{equation}\label{goal3}
|\{ \sup_{k \geq k(t,l) } \big| \sum_{i=k(t,l) }^k \sum_{I \in \mathcal{Q}_t(i)} T_I b_I^{l,r} \big| \gtrsim t^{-C} l^{-C}  \}| \lesssim t^{-C} l^{-C} \| f\|_{L^1(\mathbb{R})}
\end{equation}
for each $t,l \geq 1, r \leq l$.
We begin the argument by establishing \eqref{goal3} in the special case when $t,r \lesssim \log_2l$.

The argument is straightforward in this instance. Estimating
\[ |\{ \sup_{k \geq k(t,l) } \big| \sum_{i=k(t,l) }^k \sum_{I \in \mathcal{Q}_t(i)} T_I b_I^{l,r} \big| \gtrsim t^{-C} l^{-C} \}| \lesssim t^C l^C \sum_{k=k(t,l)}^\infty \sum_{I \in \mathcal{Q}_t(i)} \|T_I b_I^{l,r} \|_{L^1(\mathbb{R})},\]
we majorize for each $I \in \mathcal{Q}_t$
\[  \|T_I b_I^{l,r} \|_{L^1(\mathbb{R})} \lesssim 2^{r/2 - l/2} \|b_I^{l,r} \|_{L^1(\mathbb{R})}, \]
so by the disjointness of the supports of $\{ b_I^l : I \in \mathcal{Q}_t\}$, we may sum
\[ \aligned 
t^C l^C \sum_{i=k(t,l)}^\infty \sum_{I \in \mathcal{Q}_t(i)} \|T_I b_I^{l,r} \|_{L^1(\mathbb{R})} &\leq 
t^C l^C 2^{r/2 - l/2} \sum_{i=k(t,l)}^\infty \sum_{I \in \mathcal{Q}_t(i)} \|b_I^{l,r} \|_{L^1(\mathbb{R})} \\
&\leq t^C l^C 2^{r/2 - l/2} \| f\|_{L^1(\mathbb{R})}, \endaligned \]
for the desired inequality, since $t,r \lesssim \log_2 l$. Moving forward, we may therefore assume that 
\[ \max\{ t ,r \} \gg \log_2 l.\]

\subsection{Pointwise Considerations}
We are now prepared to consider the ``intensity" function, $\sum_{I} \mathbf{1}_{CI}$.
We need only select dyadic intervals $K$ of length $|K| = 2^{k(t,l)}$, so that
\[ t^{-C} l^{-C} 2^{r+t} \lesssim \sum_{I : b_I^{l,r} \neq 0} \mathbf{1}_{CI \cap K} \lesssim t^{C} l^{C} 2^{r+t} \]
holds, by pointwise considerations on the one hand and $L^1(\mathbb{R})$ considerations on the other, namely
\begin{align}
    \sum_K |\{ K : \sum_{I : b_I^{l,r} \neq 0} \mathbf{1}_{CI \cap K} \gg t^C l^C 2^{r+t} \}| &\leq t^{-C} l^{-C} 2^{-r-t} \sum_{I : b_I^{l,r} \neq 0} |I| \\
    & \qquad \lesssim t^{-C} l^{-C} \sum_I \| b_I^{l,r} \|_{L^1(\mathbb{R})} \lesssim t^{-C} l^{-C} \|f \|_{L^1(\mathbb{R})}.
\end{align} 
Moreover, we may restrict our attention to 
\begin{equation}\label{e:selK}
|\{ K : \sup_{k \geq k(t,l) } \big| \sum_{i=k(t,l) }^k \sum_{I \in \mathcal{Q}_t(i)} T_I b_I^{l,r} \big| \gtrsim t^{-C} l^{-C} \}| \gtrsim t^{-C} l^{-C} |K|;
\end{equation}
we again use the disjoint supports of the $\{b_I^l\}$ in making this reduction.

But now for any $K$ we have the following chain of inequalities,
\[ |K| \lesssim t^C l^C \int_K \sup_{k \geq k(t,l)} \big| \sum_{i=k(t,l) }^k \sum_{I \in \mathcal{Q}_t(i)} T_I b_I^{l,r} \big|^2 \lesssim t^C l^C 2^{-r-t} |K|, \]
which completes the proof, since $\max\{ t,r \} \gg \log_2 l$. The key steps in proving the estimate are replacing the rough cut-off $\mathbf{1}_K$ with a smooth pointwise majorant, $\chi_K$, satisfying the natural derivative condition 
\[ |\partial^j \chi_K| \leq C^j |K|^{-j} \mathbf{1}_{CK}, \; \; \; C = O(1), \]
the pointwise bound 
\[ |T_I b_{I}^{l,r}| \chi_K \lesssim 2^{-r-t} \mathbf{1}_{CK},\]
and the Rademacher-Menshov Lemma (below), which allows one to estimate
\[ \int_K \sup_{k \geq k(t,l)} \big| \sum_{i=k(t,l) }^k \sum_{I \in \mathcal{Q}_t(i)} T_I b_I^{l,r} \big|^2 \lesssim
\log^2( t^C l^C 2^{r+t} ) \cdot \sup_{\epsilon_I = \pm 1} \cdot  
\int \chi_K \cdot \big| \sum_{i \geq k(t,l) } \sum_{I \in \mathcal{Q}_t(i)} \epsilon_I T_I b_I^{l,r} \big|^2. \]
Here, we used crucially that at most $t^C l^C 2^{r+t}$ many scales are involved, and orthogonality methods to estimate the final integral.

We have shown that no $K$ satisfy \eqref{e:selK}, and the proof is complete.

\begin{lemma}[Theorem 10.6 of \cite{DTT}]
Let $ (X, \mu )$ be a measure space, and $ \{ \phi _j \,:\, 1\leq j \leq N\}$ a sequence of functions 
which satisfy the Bessel type inequality below, for all sequences of coefficients $c_j \in \{ \pm 1\}$, 
\begin{equation}\label{e:bessel}
\Bigl\lVert \sum_{j=1} ^{N}  c_j \phi _j \Bigr\rVert _{L ^2 (X)} \leq A .  
\end{equation}
Then, there holds 
\begin{equation}\label{e:RM}
\Bigl\lVert\sup _{1 \leq n \leq N} 
\Bigl\lvert 
 \sum_{j=1} ^{n}   \phi _j
\Bigr\rvert
\Bigr\rVert _{L ^2 (X)} \lesssim A   \log(2+ N) .  
\end{equation}
\end{lemma}

\typeout{get arXiv to do 4 passes: Label(s) may have changed. Rerun}

\end{document}